\documentclass[letter,preprint]{elsarticle}

\usepackage{amsmath,amsfonts,amssymb,amsthm}
\usepackage{mathrsfs}  
\usepackage{algorithm}
\usepackage{array}
\usepackage[caption=false,font=normalsize,labelfont=sf,textfont=sf]{subfig}
\usepackage{textcomp}
\usepackage{stfloats}
\usepackage{url}
\usepackage{verbatim}
\usepackage{graphicx}

\usepackage{microtype}
\usepackage{booktabs}
\usepackage{tabularx}
\usepackage{hyperref}

\usepackage{mathtools}
\usepackage[capitalize,noabbrev]{cleveref}
\usepackage[noend]{algpseudocode}
\usepackage{epstopdf}
\usepackage{multirow}
\usepackage{bm}
\usepackage{color}
\usepackage{wrapfig}
\usepackage[textsize=tiny]{todonotes}
\usepackage{enumitem}

\newcommand{\Stwo}{\mathbb{S}^{2}}
\newcommand{\dw}{\mathrm{d}\omega_d}
\newcommand{\dww}{\mathrm{d}\omega}
\newcommand{\Kop}{\mathcal{K}}
\newcommand{\LK}{\mathcal{L}_n\mathcal{K}}
\newcommand{\Lnf}{\mathcal{L}_n f}
\newcommand{\Pl}{\mathbb{P}_n}

\newtheorem{lemma}{Lemma}

\newtheorem{definition}{Definition}

\newcommand{\setsep}{\,\big|\,}

\begin{document}

\begin{frontmatter}
\title{A Survey on Spherical Designs: Existence, Numerical Constructions and Applications}
\author[1]{Congpei An}
\ead{andbachcp@gmail.com}

\author[2]{Xiaosheng Zhuang}
\ead{xzhuang7@cityu.edu.hk}

\affiliation[1]{organization = {School of Mathematics, Guangxi University},
			city={Nanning 530004},
			country={China}
}
\affiliation[2]{organization = {Department of Mathematics, City University of Hong Kong},
			city={Hong Kong SAR},
			country={China}
}

\begin{abstract}
This survey reviews spherical designs from a numerical-analysis perspective while retaining their connections with algebraic combinatorics and discrete geometry. A finite set $X_N\subset\mathbb{S}^d$ is a spherical $t$-design when equal-weight integration over $X_N$ is exact for all spherical polynomials of degree at most $t$. We organize the subject around three themes: cardinality bounds, existence and nonexistence; certified and large-scale numerical construction; and approximation and computational applications. Particular attention is paid to the relations among polynomial exactness, covering and separation geometry, conditioning, and numerical stability. We also discuss applications to numerical integration, interpolation, hyperinterpolation, differential and integral equations, spherical signal and image processing.
\end{abstract}


\end{frontmatter}

\section{Introduction}\label{sec:intro}
Numerical integration approximates an integral by a weighted sum of function values. The term \emph{quadrature} is usually used in one dimension, whereas \emph{cubature} refers to multidimensional integration \cite{ueberhuber1997numerical}. Classical examples include Newton--Cotes formulas \cite{cotes1722harmonia} and Gaussian quadrature \cite{gauss1814methodus}; broader accounts of quadrature and cubature may be found in \cite{milne1949numerical,hildebrand1956introduction,krylov1966handbook,stroud1971approximate,sobolev1974introduction,sobolev2006theory,engesl1980numerical,press1992numerical,davis2007methods}. Equal-weight formulas form a particularly important subclass because they avoid weight imbalance and can reduce the influence of random errors in sampled function values.

In 1874, Chebyshev \cite{chebyshev1874} studied equal-weight quadrature rules, seeking nodes $x_1,\ldots, x_t\in[-1,1]$ such that
\begin{equation}\label{eq:chebyqua}
\frac12\int_{-1}^1 p(x)\,\mathrm{d}x = \frac1t\sum_{i=1}^t p(x_i)
\end{equation}
holds for every polynomial $p$ of degree at most $t$. He computed these rules for $t=2,3,4,5,6,7$. In 1880, Radau \cite{radau1880} added the case $t=9$ and noted that the case $t=8$ involves complex nodes. More than fifty years later, Bernstein \cite{bernstein1937formules,bernstein1938systeme} proved that the equal-weight quadrature rules found by Chebyshev and Radau are the only such rules with all nodes real. Moreover, he proved that $N>t^2/16$ if one allows more real nodes $x_1,\ldots, x_N$ than the polynomial degree $t$. Krylov \cite{krylov1957proof} gave a simplified version of Bernstein's proof. Subsequent contributions were made by Kuzmin \cite{kuzmin1938distribution}, Geronimus \cite{geronimus1944gauss}, Salzer \cite{salzer1947tables}, Ullman \cite{ullman1962tchebycheff}, Kahaner \cite{kahaner1969equal}, Gautschi and Yanagiwara \cite{gautschi1974chebyshev}, Wagner and Volkmann \cite{wagner1991averaging}, and many others. A detailed survey is given by Gautschi in \cite{gautschi1976advances}. 
Equal-weight quadrature sums minimize the effect of random errors in the sampled function values $p(x_i)$, a useful feature when these errors are substantially larger than the truncation error \cite{gautschi1976advances}.

The passage from equal-weight quadrature on an interval to equal-weight cubature on the sphere has attracted sustained attention \cite{delsarte1977spherical,Yau1994,krylov1957proof,krylov1966handbook,bernstein1937quadrature}. Sobolev and his school studied these formulas from the viewpoint of numerical approximation and invariant cubature; these works form an important precursor to the modern theory of spherical designs \cite{sobolev1962cubature,sobolev1974introduction,SobolevVa97,sobolev2006theory}.

\begin{figure}[htbp!]
    \centering
    \includegraphics[width=0.45\textwidth]{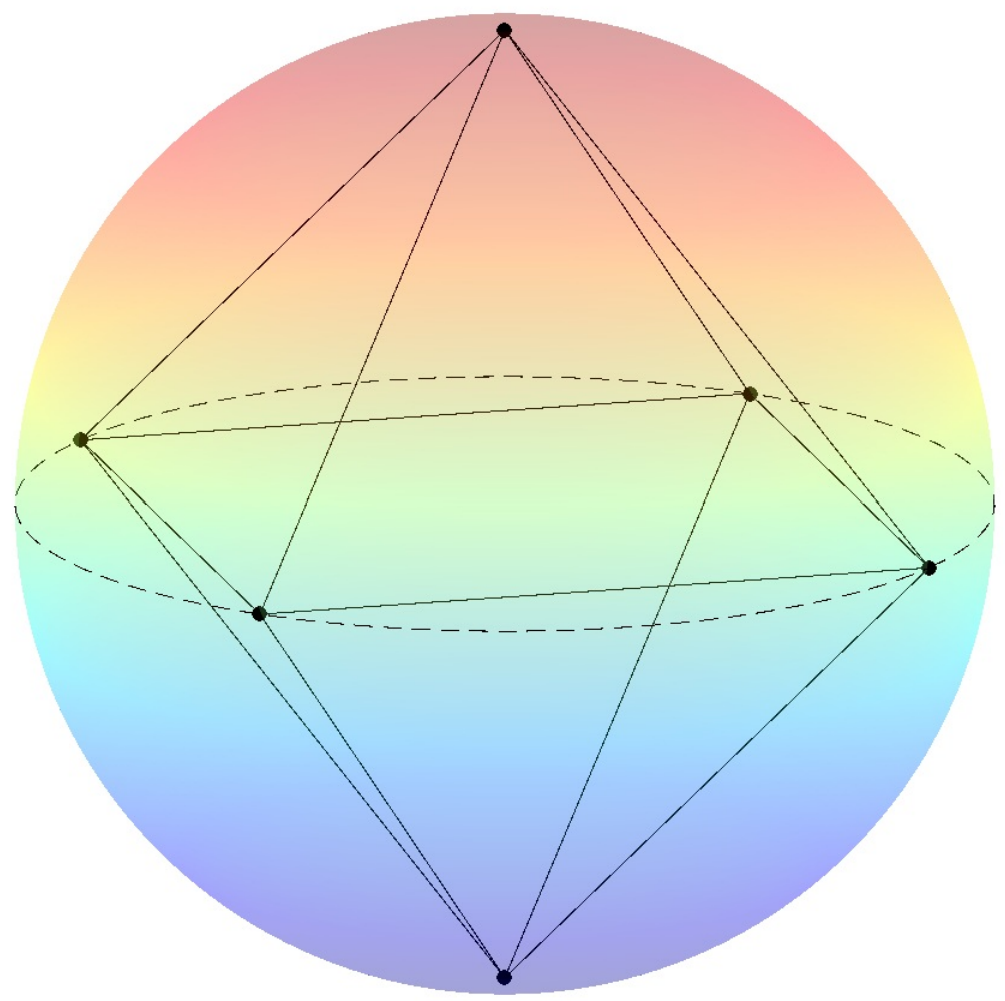}
    \includegraphics[width=0.45\textwidth]{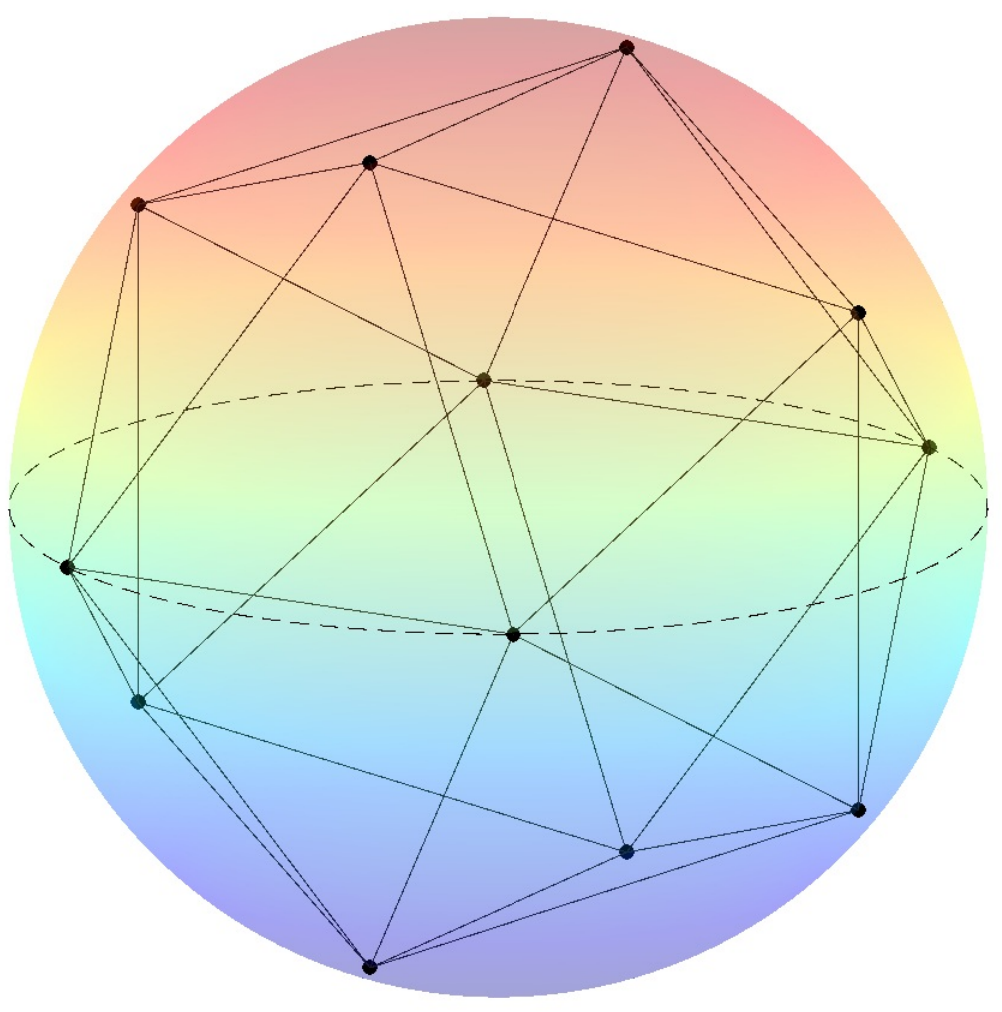}  
    \caption{Left: Spherical 3-design (cross-polytope). Right: Spherical 5-design (regular icosahedron).}
    \label{fig:std-examples}
\end{figure}

The formal concept of {\em spherical designs} was introduced by Delsarte, Goethals, and Seidel \cite{delsarte1977spherical} in 1977. Denote the unit sphere by $\mathbb{S}^d:=\{\bm x\in \mathbb{R}^{d+1}\setsep \|\bm x\|_2 =1\}$.
Let $X_{N}:=\{\bm x_1, \ldots,\bm x_N\}\subset\mathbb{S}^d$ and let $\mathbb{P}_t:= \mathbb P_t(\mathbb{S}^d)$ be the linear space of restrictions of polynomials of degree at most $t$ in $d+1$ variables to $\mathbb{S}^d$.
\begin{definition}
A point set $X_{N}\subset\mathbb{S}^d$ is a {\em spherical $t$-design} if 
\begin{equation}\label{eq:std}
    \frac{1}{|\mathbb{S}^d|}\int_{\mathbb{S}^d} p(\bm x)\,\dw(\bm x)=\frac{1}{N}\sum_{i=1}^N p(\bm x_i)\quad\quad\forall\,
    p\in \mathbb P_t,
\end{equation}
where $\dw$ is the usual (Lebesgue) surface measure on the unit sphere.
\end{definition}
The surface area of \( \mathbb{S}^d \) is denoted by \( \omega_d \),
\[
\omega_d: = |\mathbb{S}^d | = \frac{2 \pi^{\frac{d+1}{2}}}{\Gamma\left(\frac{d+1}{2}\right)}.
\]
A spherical $t$-design is therefore an equal-weight cubature rule that is exact on $\mathbb P_t$. Figure~\ref{fig:std-examples} shows a spherical $3$-design (the six vertices of a cross-polytope) and a spherical $5$-design (the twelve vertices of a regular icosahedron). The subject interacts with algebraic combinatorics \cite{bajnok1998constructions,bannai2009survey}, discrete geometry and coding theory \cite{cameron1975graph,ericson2001codes,cohn2007universally}, spherical packing \cite{musin2008kissing}, approximation theory \cite{sloan1995polynomial}, harmonic analysis \cite{fang2020theory,wang2020tight,xiao2023spherical}, machine learning \cite{lin2024sketching}, and quantum information \cite{zhu2017spherical,munemasa2006spherical}.

This survey focuses on the existence, numerical construction, and applications of spherical designs. Each section summarizes recent developments and open questions and provides references for further reading. Section 2 discusses cardinality bounds and existence, necessary conditions and nonexistence results, computer-assisted verification, and geometric properties, including covering radius and separation. Section 3 covers numerical construction methods, including optimization and fast algorithms. Section 4 explores extensions such as extremal, $t_\varepsilon$-, and nested designs. Section 5 details applications in integration, interpolation, hyperinterpolation, PDEs, integral equations, and signal processing. The paper concludes with final remarks and an extensive reference list, highlighting both theoretical advances and practical implementations in the field.

\section{Existence, nonexistence, and geometry of spherical designs}\label{sec:existence}
In this section, we review cardinality bounds and existence theorems, necessary conditions and nonexistence results, computer-assisted verification, and geometric properties of spherical designs.

\subsection{Cardinality bounds and existence}\label{subsec:Ndt}

A fundamental problem in the theory of spherical $t$-designs is to determine their minimum cardinality. We write
\[
N_{\min}(d,t):=\min\bigl\{|X|:X\subset\mathbb{S}^d \text{ is a spherical }t\text{-design}\bigr\}.
\]
The Delsarte--Goethals--Seidel lower bound \cite{delsarte1977spherical} is
\begin{equation}\label{ineq:lowerb}
    N_{\min}(d,t) \ge D(d,t) :=
    \begin{cases}
	\binom{d+k}{d}+\binom{d+k-1}{d} & \text{if } t=2k,\\[1.2ex]
	2\binom{d+k}{d} & \text{if } t=2k+1.
    \end{cases}
\end{equation}
A spherical $t$-design with $N=D(d,t)$ is called a \emph{tight} design. Bannai initiated a series of studies on the nonexistence of tight spherical designs in \cite{bannai1979tight}, establishing nonexistence results for tight $t$-designs in large dimensions. Later, in \cite{bannai1979tightI,bannai1980tightII,bannai1981uniqueness}, Bannai, Damerell, and Sloane completed the nonexistence classification of $(t,d)$ for almost all cases: If a tight $t$-design exists on $\mathbb{S}^d$ for $d\ge 2$, then $t$ must belong to $\{1,2,3,4,5,7,11\}$; moreover, if $t=11$, then $d=23$ (the $196560=D(23,11)$ minimal vectors of the Leech lattice in $\mathbb{R}^{24}$). In contrast, the 120 vertices of the regular 600-cell on $\mathbb{S}^3$ form an 11-design, but it is not tight since $D(3,11)=112$. 

On $\mathbb S^2$, the six vertices of the cross-polytope form a tight $3$-design, since $D(2,3)=6$, and the twelve vertices of the regular icosahedron form a tight $5$-design, since $D(2,5)=12$. On $\mathbb S^3$, the 24 vertices of the regular 24-cell, corresponding to the $D_4$ root system, form a spherical $5$-design but are not tight because $D(3,5)=20$. These vertices consist of the eight coordinate points $\{\pm \bm e_j:1\le j\le4\}$ and the sixteen points $\frac12(\pm1,\pm1,\pm1,\pm1)$.

A complementary problem is to find quantitative conditions on $N$, $d$, and $t$ that ensure existence for all sufficiently large $N$. In 1984, Seymour and Zaslavsky \cite{seymour1984averaging} proved that, for $N$ sufficiently large, spherical $t$-designs on $\mathbb{S}^d$ exist for all $d,t\in\mathbb{N}$. However, this proof is nonconstructive and provides no insight into how to construct practical spherical $t$-designs, nor does it establish a quantitative bound on $N$ in terms of $d$ and $t$.

Wagner and Volkmann \cite{wagner1991averaging} gave the first quantitative sufficient condition: there exists a constant $c_d$ such that, whenever $N\ge c_dt^{12d^4}$, an $N$-point spherical $t$-design on $\mathbb S^d$ exists. Here and throughout, $C_d$ and $c_d$ denote positive constants depending only on $d$. Bajnok \cite{bajnok1992construction} proved the upper bound
\[
N_{\min}(d,t)\le C_dt^{\mathcal O(d^3)},
\]
which was subsequently improved to
\[
N_{\min}(d,t)\le C_dt^{\mathcal O(d^2)}
\]
in \cite{rabau1991bounds}; see also related work of Korevaar and Meyers \cite{korevaar1993spherical,korevaar1994chebyshev} and Kuijlaars \cite{kuijlaars1993minimal}. 

Korevaar and Meyers in \cite{korevaar1993spherical} conjectured that
\[
N_{\min}(d,t) \le C_d t^d.
\]
This conjecture attracted significant interest. For instance, Kuijlaars and Saff \cite{kuijlaars1998asymptotics} emphasized its importance for $d=2$ and revealed its relation to minimal energy problems. Mhaskar, Narcowich, and Ward \cite{mhaskar2001spherical} constructed positive quadrature formulas on $\mathbb{S}^d$ using $C_d t^d$ points with nearly equal weights.

The conjecture was eventually solved by Bondarenko, Radchenko, and Viazovska through a series of works. In \cite{bondarenko2008new}, Bondarenko and Viazovska proved that
\[
N_{\min}(d,t) \le C_d t^{a_d},
\]
where the sequence $\{a_d: d=1,2,\ldots\}$ satisfies $a_1=1$, $a_2=3$, $a_{2d-1}=2a_{d-1}+d$, and $a_{2d}=a_{d-1}+a_d+d+1$. Moreover, they showed that $a_d<\frac{d}{2}\log_2(2d)$ for $d>10$. Later, in \cite{bondarenko2010spherical}, they improved the result to
\[
N_{\min}(d,t) \le C_d t^{\frac{2d(d+1)}{d+2}}
\]
based on the Brouwer fixed point theorem \cite{nirenberg1974topics}, the Marcinkiewicz--Zygmund inequality on the sphere \cite{mhaskar2001spherical}, and area-regular partitions \cite{bourgain1988distribution, kuijlaars1998asymptotics}. Finally, in \cite{bondarenko2013optimal}, using topological degree theory \cite{cho2006topological}, Bondarenko, Radchenko, and Viazovska proved the optimal bound
\[
N_{\min}(d,t) \le C_d t^d.
\]
Viazovska was awarded the Fields Medal in 2022 for her remarkable work on sphere packing and spherical designs. Dai and Feng \cite{dai2019chebyshev} extended their result to a general setting for spheres, balls, and simplexes with doubling weights.

\subsection{Necessary conditions and nonexistence results}\label{subsec:nonexistence}

The nonexistence theory for tight designs concerns the equality case in the Delsarte--Goethals--Seidel bound. However, a cardinality strictly larger than $D(d,t)$ need not be attainable. A complementary line of research, therefore, studies necessary conditions for prescribed triples $(d,t,N)$ and seeks to identify forbidden cardinalities beyond the tight case.

Boyvalenkov, Danev, and Nikova \cite{boyvalenkov1999nonexistence} derived necessary conditions for spherical designs of odd strength and odd cardinality by applying polynomial techniques to the inner-product distribution of the points. Their results yield many nonexistence statements, asymptotic improvements over the Delsarte--Goethals--Seidel estimate in the odd-cardinality setting, and detailed conclusions for strengths $3$ and $5$. Boumova, Boyvalenkov, and Danev \cite{boumova1999necessary} developed related necessary conditions in the broader setting of polynomial metric spaces, especially for designs whose cardinalities are close to the classical lower bounds.

These ideas were subsequently refined by combining polynomial inequalities for extreme inner products with additional geometric arguments. The resulting developments include new structural bounds and nonexistence results for spherical designs \cite{boumova2009polynomial}, a method tailored to odd strength and odd cardinality with applications to spherical $7$-designs \cite{boumova2009method}, improved asymptotic lower bounds for the minimum possible odd cardinality at fixed odd strength \cite{boyvalenkov2010asymptotic}, and further nonexistence results in small dimensions and asymptotic regimes \cite{boyvalenkov2013new}. This finite-parameter nonexistence theory complements the asymptotic existence results above: existence theorems identify cardinalities that are sufficient, whereas necessary conditions rule out specific parameter sets near the lower-bound regime.

The main point is that these results impose restrictions beyond the universal cardinality bound \eqref{ineq:lowerb}. Their polynomial criteria encode information about the inner-product distribution and, in the odd-strength setting, exploit parity and positivity constraints that are invisible in the Delsarte--Goethals--Seidel estimate alone. The precise auxiliary polynomials and admissible parameter ranges depend on $(d,t,N)$; the cited papers contain explicit criteria and tables of excluded parameter sets. Thus, numerical failure to locate a design at a prescribed cardinality can be compared with an independent analytic obstruction rather than interpreted as evidence of nonexistence by itself.

\subsection{Computer-assisted verification of existence}\label{subsec:interval}

Interval methods \cite{Rump2010} provide a framework for proving the existence of exact spherical $t$-designs near numerically computed configurations \cite{chen2006existence, chen2009numerical, chen2011computational}. The paper \cite{chen2011computational} presents a computational algorithm based on interval arithmetic. For a given $t$, the algorithm is designed to demonstrate the existence of a $t$-design with $(t+1)^2$ nodes on $\mathbb{S}^2$, and to compute narrow interval enclosures that are mathematically guaranteed to contain these nodes. Currently, there is no theoretical proof of the existence of a $t$-design with $(t+1)^2$ nodes for arbitrary $t$. This method establishes existence for $t=1,2,\ldots,100$ by interval verification. The spherical $t$-design is generally not unique; therefore, the method aims to identify a well-conditioned solution. It involves computing an interval enclosure for the zeros of a highly nonlinear system of dimension $(t+1)^2$. Several specialized techniques have been developed to improve the efficiency of interval arithmetic in this context. All computations are conducted using the MATLAB toolbox INTLAB \cite{Rump2010}.

In \cite{an2010well}, well-conditioned spherical $t$-designs on $\mathbb S^2$ are constructed and verified using interval methods up to $t=60$. These methods also provide a guaranteed interval for the determinant, enhancing the reliability of the results and confirming the geometric properties of the resulting spherical designs, such as separation and mesh norm. Using the supercomputing center at Jinan University, well-conditioned spherical $t$-designs have been constructed and verified using interval methods up to $t=160$ with $N=(t+1)^2$, as detailed in \cite{Chensy2017} and \cite{anchen2016Numerical}. In an invited lecture at the Forum for Women Scientists, ICCM2025, Shanghai, Chen \cite{Chen_2026} emphasized the importance and feasibility of interval methods for proving the existence of spherical $t$-designs.

\subsection{Separation and covering geometry}\label{subsec:geom}

Throughout this subsection, we use the geodesic distance on $\mathbb{S}^d$ defined by
\[
\operatorname{dist}(\bm x,\bm y) := \cos^{-1}(\bm x\cdot \bm y), \qquad \bm x,\bm y\in \mathbb{S}^d,
\]
where $\bm x\cdot \bm y$ denotes the standard Euclidean inner product in $\mathbb{R}^{d+1}$.

Given a finite set of distinct points $X_N=\{\bm x_1,\ldots,\bm x_N\}\subset \mathbb{S}^d$, the quality of the point distribution is characterized by the following quantities.

\begin{definition}
The mesh norm $h_{X_N}$ of a point set $X_N\subset\mathbb{S}^d$ is
\begin{equation}
    h_{X_N} := \max_{\bm y\in\mathbb{S}^d} \min_{\bm x_i \in X_N} \cos^{-1}(\bm y\cdot \bm x_i).
\end{equation}
\end{definition}

In other words, the mesh norm can be regarded as the geodesic radius of the largest hole in the mesh \cite{an2010well}; equivalently, it is the covering radius, i.e., the smallest radius for which $N$ identical spherical caps centered at the points $\bm x_i$ cover the sphere \cite{sole1991covering}.

For comparison with the linear-programming literature, it is useful to introduce the equivalent inner-product covering parameter
\begin{equation}\label{eq:inner-product-covering}
    \eta(X_N):=\min_{\bm y\in\mathbb{S}^d}\max_{1\le i\le N}\bm y\cdot\bm x_i.
\end{equation}
Since $\cos^{-1}$ is decreasing on $[-1,1]$, one has
\[
    h_{X_N}=\cos^{-1}\!\bigl(\eta(X_N)\bigr).
\]
Fazekas and Levenshtein \cite{fazekas1995upper} obtained linear-programming bounds for the corresponding covering parameter in polynomial metric spaces; in the spherical setting their classical bound depends on the dimension and the strength, but not on the cardinality. Boyvalenkov and Stoyanova \cite{boyvalenkov2021covering} subsequently derived cardinality-sensitive lower and upper bounds for $\eta(X_N)$. Schematically, their estimates take the form
\[
   \eta^-_{\mathrm{LP}}(d,t,N)\le \eta(X_N)\le \eta^+_{\mathrm{LP}}(d,t,N),
\]
where the two endpoints are obtained from signed measures and associated systems of orthogonal polynomials; their explicit forms and admissible parameter ranges depend on the parity of $t$. For designs of even strength, the lower estimate can improve the Fazekas--Levenshtein bound, while the upper-bound method also applies in several odd-strength and antipodal settings. In the geodesic convention used here, a lower bound for $\eta(X_N)$ is an upper bound for $h_{X_N}$, whereas an upper bound for $\eta(X_N)$ is a lower bound for $h_{X_N}$.

\begin{definition}
The separation distance of a point set $X_N\subset \mathbb{S}^d$ is defined by
\begin{equation}
    \delta_{X_N} := \min_{\bm x_i,\bm x_j\in X_N,\, i\neq j} \cos^{-1}(\bm x_i\cdot \bm x_j).
\end{equation}
\end{definition}
Note that the separation distance can also be defined using the Euclidean distance $|\bm x -\bm y|$, since $|\bm x-\bm y|\le \cos^{-1}(\bm x\cdot \bm y) \le\frac{\pi}{2}|\bm x -\bm y|$ for $\bm x,\bm y\in \mathbb S^d$. 

A spherical $t$-design $X_N$ is said to be \emph{well-separated} if
\[
 \delta_{X_N} \ge \lambda_d N^{-1/d},
\]
for some constant $\lambda_d$. This condition is usually considered in the regime $N\ge t^d$ (or $N$ of order $t^d$). If two points $\bm x_i$ and $\bm x_j$ are very close, the condition number of the interpolation matrix becomes large. This is why well-separated spherical $t$-designs are desirable. Numerical evidence for the existence of well-separated spherical $t$-designs can be found in \cite{an2010well}. It was shown in \cite{bondarenko2015well} that there exist well-separated spherical $t$-designs with $N$ points of the optimal order $t^d$.

The mesh ratio $\rho_{X_N}$ defined by
\[
\rho_{X_N} := \frac{2h_{X_N}}{\delta_{X_N}} \ge 1
\]
is a good measure of the quality of the geometric distribution of $X_N$: the smaller $\rho_{X_N}$ is, the more uniformly the points are distributed on $\mathbb{S}^d$ \cite{leopardi2007distributing, an2010well}. Equivalently, $\rho_{X_N}$ can be interpreted as the ratio of the covering radius $h_{X_N}$ to the packing radius $\delta_{X_N}/2$.

\section{Numerical construction of spherical designs}\label{sec:numerical}
In this section, we focus on the numerical construction of spherical designs on $\mathbb{S}^d$. 

\subsection{Optimization problems for finding spherical designs}
To find spherical $t$-designs for a given $t$, one may formulate the point set $X_N$ as the solution to a system of nonlinear equations. This naturally leads to various optimization formulations, which in turn call for the development of specialized numerical algorithms.

We denote by
\[
\{Y_{\ell}^{k}: k=1,\ldots,Z(d,\ell),\; \ell=0,1,2,\ldots\}
\]
an $L^2$-orthonormal system on $\mathbb S^d$ of real spherical harmonics of degree $\ell$ and order $k$, normalized with respect to the surface measure $\mathrm{d}\omega_d$. Here $Z(d,\ell)$ is the dimension of the space of spherical harmonics of degree $\ell$. Thus,
\[
\int_{\mathbb S^d}Y_\ell^k(\bm x)Y_{\ell'}^{k'}(\bm x)\,\mathrm{d}\omega_d(\bm x)
=\delta_{\ell\ell'}\delta_{kk'}.
\]
For $d=2$, one has $Z(2,\ell)=2\ell+1$.

\begin{lemma}[\cite{delsarte1977spherical, sloan2009variational, Cui1997Eq}]
A finite set $X_N=\{\bm x_1,\ldots,\bm x_N\}\subset \mathbb{S}^d$ is a spherical $t$-design if and only if the so-called Weyl sums vanish:
\begin{equation} \label{lemma:eqcon}
    \sum_{j=1}^{N} Y_{\ell}^{k}(\bm x_j)=0,\qquad k=1,\ldots,Z(d,\ell),\quad \ell=1,\ldots,t.
\end{equation}
\end{lemma}

One way to convert this system into an optimization problem is to consider a nonnegative sum of squares of the Weyl sums. Sloan and Womersley \cite{sloan2009variational} introduced the variational characterization
\begin{equation}\label{eq:A}
   A_{N,t}(X_N):=\frac{|\mathbb S^d|}{N^2}\sum_{\ell=1}^t\sum_{k=1}^{Z(d,\ell)}
   \left|\sum_{j=1}^N Y_{\ell}^{k}(\bm x_j)\right|^2.
\end{equation}
It was shown in \cite{sloan2009variational} that a point set $X_N\subset \mathbb{S}^d$ is a spherical $t$-design if and only if $A_{N,t}(X_N)=0$. Since $A_{N,t}\ge 0$, any $t$-design is a global minimizer of $A_{N,t}$ with value $0$. Conversely, if the global minimum of $A_{N,t}$ is strictly positive, then no $N$-point spherical $t$-design exists. Thus the problem of seeking $X_N$ can be transformed into the optimization problem of minimizing $A_{N,t}$:
\[
X_N = \arg\min_{X_N\subset \mathbb{S}^d} A_{N,t}(X_N).
\]
Here $A_{N,t}$ is a function on the product manifold $(\mathbb{S}^d)^N$.

It is important to distinguish the preceding mathematical equivalence from what a numerical computation can certify. The failure of an optimization algorithm to reduce $A_{N,t}$ to zero does not by itself prove nonexistence, since the algorithm may have converged to a nonglobal local minimum or may have stalled before reaching a zero. A rigorous nonexistence conclusion requires either a certified positive lower bound for the global minimum or analytic necessary conditions such as those reviewed in \cref{subsec:nonexistence}.

\subsection{Newton-like method}
For a general optimization problem, it is typically difficult to guarantee that a numerical algorithm will converge to a global minimizer. In this regard, Sloan and Womersley \cite{sloan2009variational} examined the first-order optimality conditions of $A_{N,t}$, and showed that stationary configurations satisfying a condition on the mesh norm are spherical designs and hence global minimizers. The main result in \cite{sloan2009variational} states that if $X_N\subset \mathbb{S}^d$ is a stationary point set of $A_{N,t}$ and its mesh norm satisfies
\[
h_{X_N} < \frac{1}{t+1},
\]
then $X_N$ is a spherical $t$-design. This result opened a pathway to proving the existence of spherical $t$-designs with a number $N$ of points proportional to $(t+1)^d$ for fixed $d$. Using a quasi-Newton method, a numerical example with $d=2$ and $t=19$ suggested that computational minimization of $A_{N,t}$ could be a valuable tool for discovering new spherical designs for moderate and large $t$.

In 2018, Womersley \cite{womersley2018efficient} considered the generation of efficient spherical $t$-designs with good geometric properties, where ``efficient'' refers to the number of points $N$ being comparable to $\dim\mathbb P_t(\mathbb S^d)$. For fixed $d$ and $t\to\infty$,
\[
\dim\mathbb P_t(\mathbb S^d)
=\binom{t+d}{d}+\binom{t+d-1}{d}
\sim\frac{2}{d!}t^d.
\]
Good geometric properties are characterized by a low mesh ratio $\rho_{X_N}$. The results for $\mathbb{S}^2$ include computed spherical $t$-designs for $t=1,\ldots,180$ and symmetric (antipodal) $t$-designs for degrees up to 325, all exhibiting low mesh ratios. These point sets provide excellent candidates for numerical integration on the sphere. Furthermore, the methods can also be used to explore spherical $t$-designs computationally for $d=3$ and higher. Such work is based on quasi-Newton methods and utilizes supercomputer systems.

\subsection{Barzilai--Borwein method}\label{subsec:BB}

Newton-like methods use Hessian information or approximations to it, which may be computationally expensive for large-scale problems. Newton steps are typically obtained by solving a linear system, without explicitly forming the inverse of the Hessian. The Barzilai--Borwein method \cite{YuanOpt2006} is a gradient-based method with modified step sizes, motivated by Newton's method but without involving the Hessian. In \cite{an2020numerical}, it is shown that if $X_N$ is a stationary point set of $A_{N,t+1}$ and the smallest singular value of the basis matrix (i.e., the spherical harmonic evaluation matrix $\bm Y_{t+1}(X_N)$, see Definition~\ref{defH}) is positive, then $X_N$ is indeed a spherical $t$-design. Numerical examples are presented for spherical $t$-designs with $N=(t+2)^2$ points for values of $t+1$ up to 127 on $\mathbb S^2$. As the authors point out, this method demonstrates high efficiency and accuracy. Under the stated assumptions, exact stationarity and positivity of the smallest singular value imply that the configuration is a spherical $t$-design and hence a global minimizer of $A_{N,t}$. In numerical computations, stationarity residuals and computed singular values provide numerical checks; rigorous certification requires validated bounds for the relevant hypotheses. Numerical experiments indicate that the Barzilai--Borwein method is faster than the quasi-Newton method when constructing spherical $t$-designs with large numbers of nodes.

An and Xu \cite{AnXu2025} propose a regularized Barzilai--Borwein (RBB) method, which tunes the spectral gradient step sizes by constructing a regularized least-squares model. They employ the RBB method to numerically compute spherical $t$-designs up to $t=130$ with $N=(t+1)^2$ on $\mathbb S^2$. Xu \cite{Xu2025} introduces an interpolated least-squares model, which unifies the two classical BB step sizes by modifying the degree of the independent variable in the original least-squares model. The results of \cite{AnXu2025} and \cite{Xu2025} provide new perspectives for the study of BB-like methods, and by appropriately selecting parameters, both approaches can deliver efficient step sizes for gradient descent algorithms.

\subsection{Newton and conjugate gradient methods on manifolds}\label{subsec:newton}

The sphere $\mathbb{S}^d$ is a compact Riemannian manifold. Consequently, numerical optimization methods on manifolds can be applied to find spherical $t$-designs. In \cite{graf2011computation}, the authors consider the problem of finding numerical spherical $t$-designs on $\mathbb{S}^2$ for high polynomial degrees $t\in\mathbb{N}$. Specifically, they compute local minimizers of the quadrature error $A_{N,t}$. The same variational characterization was used by Sloan and Womersley in \cite{sloan2009variational}. For the minimization problem, the authors employ several nonlinear optimization methods on manifolds, such as Newton's method and conjugate gradient methods.

The paper \cite{graf2011computation} demonstrates that, by utilizing the Nonequispaced Fast Spherical Fourier Transform (NFSFT) on $\mathbb S^2$, gradient and Hessian evaluations can be performed in $O(t^2\log t + N\log^2(1/\epsilon))$ arithmetic operations, where $\epsilon>0$ is a prescribed accuracy. Using these methods, the authors present numerical spherical $t$-designs on $\mathbb S^2$ for $t\le 1000$, even in cases when $N\approx t^2/2$.

\subsection{Trust-region methods}\label{subsec:trust}

Trust-region methods \cite{YuanOpt2006} are a class of optimization techniques. Unlike line-search methods, which perform a line search in each iteration, trust-region methods compute a trial step by solving a subproblem that minimizes a model function within a trust region. This approach allows for the use of nonconvex models, making these methods suitable for nonconvex and ill-conditioned problems.

The function $A_{N,t}$ is generally nonconvex. In \cite{xiao2023spherical}, the authors numerically solve the minimization problem using the trust-region method. The function $A_{N,t}$ is formulated as a band-limited expression involving the spherical harmonic evaluation matrix $\bm Y_t:=\bm Y_t(X_N)$ and its adjoint $\bm Y_t^\star$, allowing compact matrix-vector forms for $A_{N,t}$, its gradient $\nabla A_{N,t}$, and its Hessian $\mathcal H(A_{N,t})$. To accelerate computations, NFSFTs are used for matrix-vector operations involving $\bm Y_t$ and $\bm Y_t^\star$. Using FFT-based spherical harmonic transforms reduces the computational cost and makes the method feasible for large $N\approx (t+1)^2$ on $\mathbb S^2$. This efficiency is crucial in the trust-region optimization framework, where $A_{N,t}$ and its derivatives must be evaluated repeatedly. Through the fast spherical harmonic transform, each iteration's cost can be kept to $\mathcal O(t^2\log^2 t + N\log^2(1/\epsilon))$, thus facilitating efficient large-scale computations on $\mathbb S^2$. The paper \cite{xiao2023spherical} finds numerical spherical $t$-designs on $\mathbb S^2$ with $t$ as large as 3000.

\subsection{Comparison of numerical construction methods}

Table~\ref{tab:comparison} summarizes the main computational features of the numerical construction methods discussed above. Here $Z_t = \dim(\mathbb P_t(\mathbb S^2)) = (t+1)^2$ denotes the dimension of the polynomial space.

\begin{table}[htbp]
\centering
\caption{Comparison of numerical construction methods for spherical $t$-designs on $\mathbb S^2$.}
\label{tab:comparison}
\begin{scriptsize}
\begin{tabularx}{\textwidth}{@{}l l l X@{}}
\toprule
Method & Per-iteration cost & Storage & Typical local behavior \\
\midrule
Newton-like  & $O(N Z_t)$ & $O(Z_t^2)$ (Hessian) & Superlinear under standard assumptions \\
Barzilai--Borwein & $O(N Z_t)$ & $O(N)$ & Usually linear \\
Manifold Newton/CG & $O(t^2\log t + N\log^2(1/\epsilon))$ & $O(N)$--$O(N^2)$ & Newton: quadratic; CG: problem-dependent \\
Trust-region & $O(t^2\log^2 t + N\log^2(1/\epsilon))$ & $O(N Z_t)$ & Problem-dependent; often superlinear locally \\
\bottomrule
\end{tabularx}
\end{scriptsize}
\end{table}

The choice of method involves a trade-off between available derivative information, storage requirements, and convergence speed. The complexity and convergence entries in Table~\ref{tab:comparison} are implementation-dependent summaries of the behavior reported in the cited numerical studies; they should not be read as unconditional guarantees. Limited-memory quasi-Newton methods with low-rank Hessian approximations are more directly comparable to the Barzilai--Borwein method in terms of storage and per-iteration cost.

\section{Miscellaneous extensions of spherical designs}\label{sec:extension}
In this section, we focus on some extensions of spherical designs, including extremal spherical designs, spherical $t_\varepsilon$-designs, nested spherical designs, and designs on local spherical regions.

\subsection{Extremal spherical designs}
The formal definition of \emph{extremal spherical designs} can be found in \cite{an2010well}. 

\begin{definition}\label{defH}
A set $X_N =\{\bm x_1,\ldots,\bm x_N\}\subset\mathbb{S}^2$ of $N\ge (t+1)^2$ points is an extremal spherical $t$-design if the determinant of the matrix
\[
H_t(X_N) := \bm Y_t(X_N)\bm Y_t(X_N)^T \in \mathbb{R}^{(t+1)^2\times (t+1)^2}
\]
is maximal subject to the constraint that $X_N$ is a spherical $t$-design, where $\bm Y_t(X_N):=(Y_\ell^k(\bm x_i))_{0\le \ell\le t,\,1\le k\le 2\ell+1;\,1\le i\le N}$ is the spherical harmonic evaluation matrix of size $(t+1)^2\times N$.
\end{definition}

The paper \cite{an2010well} presents a method for constructing extremal spherical designs, which gives \emph{well-conditioned spherical designs} with $N\ge (t+1)^2$ points on $\mathbb S^2$. A crucial aspect of this approach is the application of interval methods, which are vital for demonstrating the existence of a true spherical $t$-design that closely approximates the computed points.

\subsection{Spherical $t_\varepsilon$-designs}
Although interval methods have confirmed the existence of spherical $t$-designs on $\mathbb S^2$ for $t\le 100$ with $N=(t+1)^2$, constructing designs with equal positive weights and polynomial exactness of degree $t$ remains challenging in practice for large $t$. Zhou and Chen \cite{zhou2018spherical} relax the equal-weight constraint: the weights may vary within a prescribed interval, while the rule remains exact for all polynomials of degree at most $t$.

\begin{definition}
A spherical $t_\varepsilon$-design with $0\le \varepsilon<1$ on $\mathbb{S}^d$ is a set of points $X_{N,\varepsilon}=\{\bm x_{\varepsilon,1},\ldots,\bm x_{\varepsilon,N}\}\subset\mathbb{S}^d$ together with weights $w_1,\ldots,w_N$ such that the quadrature rule
\[
\int_{\mathbb{S}^d}p(\bm x)\,\dw(\bm x)=\sum_{i=1}^{N}w_i p(\bm x_{\varepsilon,i})
\]
is exact for all spherical polynomials of degree at most $t$. The weights satisfy
\[
\frac{|\mathbb{S}^d|}{N}(1-\varepsilon)
\le w_i\le
\frac{|\mathbb{S}^d|}{N}(1-\varepsilon)^{-1},
\qquad i=1,\ldots,N.
\]
\end{definition}

Since $\varepsilon<1$, the weights are positive. Moreover, exactness for constants (polynomials of degree zero) forces $\sum_{i=1}^N w_i=|\mathbb{S}^d|$, i.e., the mean of the weights is exactly $|\mathbb{S}^d|/N$.

Spherical $t_\varepsilon$-designs serve as a bridge between true spherical $t$-designs and positive-weight quadrature rules. Since every spherical $t$-design is also a spherical $t_\varepsilon$-design for any $0\le\varepsilon<1$, the existence of $t_\varepsilon$-designs follows from that of $t$-designs. The more interesting question is for which (smaller) values of $N$ a $t_\varepsilon$-design exists; the additional freedom in the weights suggests that $t_\varepsilon$-designs may exist with fewer points than spherical $t$-designs. Numerical experiments in \cite{zhou2018spherical} demonstrate that by increasing $\varepsilon$, it is possible to maintain degree-$t$ polynomial exactness with fewer points than those required by spherical $t$-designs. Related generalizations were studied in \cite{zhou2020quadrature}.

\subsection{Nested spherical designs}
On the sphere, function approximation approaches such as hyperinterpolation \cite{sloan1995polynomial}, multiscale analysis \cite{gia2010multiscale}, and localized systems \cite{narcowich2006localized, legia2009localized} rely on the polynomial exactness of spherical designs. In the approximation schemes considered here, reducing the error generally requires increasing the design strength $t$ and, consequently, the number of nodes. In real-world applications such as signal processing or super-resolution imaging on the sphere, the design strength may need to be large to match the data resolution and achieve the desired accuracy. Consequently, the total number of nodes in a spherical design becomes a major determinant of computational cost and storage requirements at high resolution. Furthermore, spherical framelets such as those in \cite{wang2020tight, xiao2023spherical} utilize a sequence of high-order cubature rules to achieve multiscale analysis and localization. In \cite[Sec. 4.5]{womersley2018efficient}, the importance of the nested property of spherical designs was discussed. Nested structures of cubature rules could significantly reduce storage and computational burdens. It is therefore desirable to have nested spherical designs for both theory and applications.

\begin{definition}
Let $X_{t_1,N_1}\subset \mathbb S^d$ and $X_{t_2,N_2}\subset \mathbb S^d$ be a spherical $t_1$-design and a spherical $t_2$-design, respectively. They are called {\em nested} if $X_{t_1,N_1}\subset X_{t_2,N_2}$ with $t_1<t_2$.
\end{definition}

Given the above discussion, it is natural to ask the following important questions.
\begin{itemize}
\item[(Q1)] Do nested spherical designs exist?
\item[(Q2)] How many nodes are needed to extend a $t_1$-design to a $t_2$-design ($t_2>t_1$)?
\item[(Q3)] What is the optimal bound for the number of points needed to extend a $t_1$-design to a $t_2$-design?
\end{itemize}
In \cite{zheng2024existence}, the authors answered these questions. They proved the existence of nested spherical designs for sufficiently large $N$, and established an upper bound of order $N=\mathcal O(t^{2d+1})$ that guarantees the existence of nested spherical designs. It is believed that such a bound is not optimal, and it is conjectured that the optimal bound should be of order $N=\mathcal O(t^d)$. For some special cases $t_1<t_2$, the paper \cite{zheng2024existence} shows that there exist nested spherical designs with cardinalities of order $c_d t_1^d$ and $c_d t_2^d$, respectively, with the same constant $c_d$.

\subsection{Designs on local spherical regions}
Positive-weight cubature formulas and equal-weight designs are distinct notions, even though both are governed by polynomial exactness and often use Marcinkiewicz--Zygmund inequalities as a principal tool \cite{mhaskar2001spherical}. For a local region $\Omega\subset\mathbb S^d$, a design must be defined relative to the restricted surface measure on $\Omega$ and the corresponding polynomial space. Recent work develops such equal-weight designs for spherical caps and spherical zones, proves existence under suitable geometric assumptions, and studies their numerical construction and approximation properties \cite{liChen2024,LiChen2025}. These results should therefore be viewed as local analogues of spherical designs, rather than as automatic consequences of the existence of positive-weight cubature formulas.

\section{Applications}

Once (numerical) spherical $t$-designs are obtained, they can be applied in various areas, e.g., numerical integration and function approximation \cite{an2010well, hesse2010numerical, an2012regularized, hesse2012numerical, xiao2025spherical}, multiscale analysis and signal processing \cite{narcowich2006localized, legia2009localized, wang2020tight, xiao2023spherical}, interpolation and hyperinterpolation \cite{sloan1995polynomial}, spherical partial differential equations \cite{wu2023restricted}, and spherical integral equations \cite{MR878693}.

\subsection{Numerical integration and function approximation}

The most direct applications of spherical $t$-designs are numerical integration \cite{an2010well, hesse2012numerical} and function approximation \cite{xiao2025spherical}. In \cite{an2010well}, an example of approximating an exponential function on the sphere is presented. In \cite{xiao2025spherical}, the approximation of smooth and non-smooth functions by spherical harmonics using spherical designs is investigated. The decay behavior of the error is consistent with the error bounds for QMC designs in the literature \cite{brauchart2014qmc}:
\begin{equation}\label{equ:stderror}
\sup_{\substack{f\in H^s(\mathbb{S}^d),\\ \|f\|_{H^s}\le 1}}
\left\lvert \frac{|\mathbb{S}^d|}{N}\sum_{j=1}^N f(\bm x_j)-\int_{\mathbb{S}^d} f(\bm x)\,\dw(\bm x) \right\rvert
\le \frac{C(s,d)}{t^s},
\end{equation}
where $H^s(\mathbb{S}^d)$ is the Sobolev space with smoothness $s>d/2$, so that point evaluation is continuous, and $C(s,d)$ is independent of the particular design. For a sequence of spherical $t$-designs satisfying the optimal cardinality relation $N\asymp t^d$, the bound $t^{-s}$ is equivalent to the QMC rate $N^{-s/d}$. Such sequences provide fundamental examples of QMC designs. The numerical examples in \cite{brauchart2014qmc} illustrate their competitive integration performance. Moreover, \cite{anchen2016Numerical} and \cite{Chensy2017} employed well-conditioned spherical $t$-designs to approximate integrals of singular functions with the help of homeomorphic transformations: the Atkinson transform and the Sidi transform (see \cite{atkinson2005quadrature, sidi2005analysis} for details). Numerical examples demonstrate the clear advantage of spherical $t$-designs compared with other quadrature rules and node configurations, including the bivariate trapezoidal rule and equal-area partition points; see \cite[Sec. 5]{anchen2016Numerical} and \cite{Chensy2017}. It should be noted, however, that spherical $t$-designs cannot be generated on the fly for arbitrarily large $t$; their construction requires dedicated numerical optimization.

\subsection{Polynomial interpolation}
For a given $f\in C(\mathbb{S}^2)$, the classical interpolation operator $\Lambda_t$ is defined by
\[
\Lambda_t f\in \mathbb{P}_t,\qquad \Lambda_t f(\bm x_j)=f(\bm x_j),\quad j=1,\ldots,N,
\]
where $N=\dim(\mathbb P_t(\mathbb S^2))$. The interpolant is expressed as
\begin{equation}
    \Lambda_t f=\sum_{j=1}^{N} f(\bm x_j)\,\ell_j,
\end{equation}
where $\ell_j\in\mathbb P_t$ are the Lagrange polynomials satisfying $\ell_j(\bm x_i)=\delta_{ij}$, $i,j=1,\ldots,N$. The Lebesgue constant for interpolation is defined by
\begin{equation}
    \label{LebesgueInterpolation}
    \|\Lambda_t\| := \sup_{f\in C(\mathbb{S}^2)\setminus\{0\}} \frac{\|\Lambda_t f\|_\infty}{\|f\|_\infty}
    = \max_{\bm x\in\mathbb{S}^2}\sum_{j=1}^{N} |\ell_j(\bm x)|.
\end{equation}

For the well-conditioned spherical $t$-designs computed in \cite{an2010well}, Figure~6.1 of that paper reports Lebesgue constants between $\sqrt{t}$ and $(t+1)^2$, with values relatively close to $t+1$. Note that $\sqrt{t}$ is the growth rate of the Lebesgue constant for the orthogonal projection \cite{an2012regularized}, which is known to be a lower bound for the Lebesgue constant for interpolation \cite{sloan2004extremal}. A nonlinear fit to the data in \cite[Figure 6.1]{an2010well} gives $0.8025(t+1)^{1.12}$ as an empirical estimate of the Lebesgue constant.

\subsection{Hyperinterpolation}
Hyperinterpolation, a discrete analogue of the $L_2$ orthogonal projection, is a powerful technique for multivariate polynomial approximation introduced by Ian H. Sloan in 1995 \cite{sloan1995polynomial}. The hyperinterpolation operator $\mathcal L_L$ is constructed by replacing the Fourier integrals in the $L_2$ orthogonal projection onto $\mathbb{P}_L$ with a quadrature rule that exactly integrates all spherical polynomials of degree up to $2L$. It is known \cite{sloan1995polynomial} that for $L\ge 3$, the number of quadrature points required for hyperinterpolation must exceed the dimension of the polynomial space. Thus, hyperinterpolation stands in contrast to interpolation, being a discrete projection method rather than a pointwise matching scheme. In what follows, we restrict the quadrature rules employed in hyperinterpolation to be spherical designs \cite{an2012regularized}.

Given a spherical $t$-design $X_N = \{\bm x_1,\ldots,\bm x_N\}\subset\mathbb{S}^2$ with $t\ge 2L$, we define the semi-inner product $(\cdot,\cdot)_N$ for two continuous functions $f,g\in C(\mathbb{S}^2)$ by
\begin{equation}\label{hyper-inner}
(f,g)_N := \frac{4\pi}{N} \sum_{j=1}^{N} f(\bm x_j)\, g(\bm x_j).
\end{equation}
It follows that for any $p,q\in\mathbb{P}_L$,
\[
(p,q)_N = (p,q)_{L_2} = \int_{\mathbb{S}^2} p(\bm x) q(\bm x)\, \mathrm{d}\omega_2(\bm x),
\]
since $pq\in\mathbb{P}_{2L}(\mathbb{S}^2)$ and $t\ge 2L$.  Note that for $f\in C(\mathbb{S}^2)$, $(f,f)_N=0$ implies $f(\bm x_j)=0$ for $j=1,\ldots,N$, but does not imply $f\equiv 0$. Thus \eqref{hyper-inner} induces only a seminorm $\|f\|_N := \sqrt{(f,f)_N}$ on $C(\mathbb{S}^2)$.

The hyperinterpolant of a function $f\in C(\mathbb{S}^2)$ is defined by
\[
\mathcal L_L f(\bm x) = \sum_{\ell=0}^{L} \sum_{k=1}^{2\ell+1} (f, Y_{\ell}^k)_N\, Y_{\ell}^k(\bm x), \qquad \bm x\in\mathbb{S}^2.
\]
In particular, $\mathcal L_L f\in\mathbb{P}_L$. Using the reproducing kernel
\begin{equation}\label{equ:reprodkernel}
G_L(\bm x,\bm y) = g_L(\bm x\cdot\bm y)
= \sum_{\ell=0}^{L} \sum_{k=1}^{2\ell+1} Y_{\ell}^k(\bm x)Y_{\ell}^k(\bm y)
= \sum_{\ell=0}^{L} \frac{2\ell+1}{4\pi} P_\ell(\bm x\cdot\bm y),
\end{equation}
the hyperinterpolant can be expressed as
\begin{align*}
\mathcal L_L f(\bm x) &= (f, G_L(\bm x,\cdot))_N
= \frac{4\pi}{N} \sum_{j=1}^{N} f(\bm x_j)\, g_L(\bm x\cdot\bm x_j) \\
&= \sum_{\ell=0}^{L} \frac{2\ell+1}{N} \sum_{j=1}^{N} f(\bm x_j)\, P_\ell(\bm x\cdot\bm x_j), \qquad \bm x\in\mathbb{S}^2,
\end{align*}
which is the discrete analogue of the orthogonal projection
\[
\mathcal P_L f(\bm x) = (f, G_L(\bm x,\cdot))_{L_2}
= \int_{\mathbb{S}^2} f(\bm y)\, g_L(\bm x\cdot\bm y)\,\mathrm{d}\omega_2(\bm y).
\]

Spherical $t$-designs have been extensively applied in hyperinterpolation and its variants \cite{anwu2024hyperinterpolation, AnRanWu2025, an2012regularized, Lin2021, lin2024sketching, an2025hybrid, an2025hard, an2021lasso, an2022quadrature}. Due to space limitations, we only mention the work on approximating singular and oscillatory functions on the sphere via hyperinterpolation using spherical $t$-designs as nodes \cite{anwu2024hyperinterpolation}. The efficient hyperinterpolation method of An and Wu \cite{anwu2024hyperinterpolation} on $\mathbb{S}^2$ uses spherical $t$-designs to define the cubature rule for approximating functions of the form
\[
F = K f,
\]
where $K$ is a kernel function that is typically singular or highly oscillatory. The key idea is to {\em precompute the modified moments}
\begin{equation}\label{equ:modifiedmoments}
    \int_{\mathbb{S}^2} K\, Y^k_{\ell}\, \mathrm{d}\omega_2
\end{equation}
analytically or via stable iterative procedures, thereby avoiding direct numerical integration of the singular or oscillatory kernel $K$ and significantly reducing the required number of quadrature points. For this variant, let $n$ denote the output degree, corresponding to $L$ above, and relax the cubature exactness to $t=n+n'$ with $0<n'<n$. Define
\[
E_m(g):=\inf_{p\in\mathbb P_m}\|g-p\|_\infty,
\]
and let $\chi^*\in\mathbb P_{n'}$ be a best uniform approximation to $f$, so that $\|f-\chi^*\|_\infty=E_{n'}(f)$. Under the assumptions of \cite{anwu2024hyperinterpolation}, for continuous $K$, the $L^2$ error bounds have the form
\begin{align*}
\|\mathcal S_n(Kf)-Kf\|_2
&\lesssim E_{n'}(f)+E_n(K\chi^*),\\
\|\mathcal L_n(Kf)-Kf\|_2
&\lesssim E_{n'}(Kf),
\end{align*}
where $\mathcal S_n$ denotes efficient hyperinterpolation and the implied constants depend on the kernel and the stated stability assumptions. These bounds distinguish the approximation of the smooth factor $f$ from that of the product $Kf$. Integrable singular kernels are treated separately in \cite{anwu2024hyperinterpolation} through $L^1$ and $L^2$ estimates.

For oscillatory spherical harmonics $Y_\ell^k$ and singular kernels such as $|\bm\xi-\bm x|^\nu$ or $\log|\bm\xi-\bm x|$, the numerical experiments in \cite{anwu2024hyperinterpolation} show improved accuracy and robustness for the tested kernels and point configurations using the same number of sampling points.

For a detailed survey on hyperinterpolation, we refer the reader to \cite{AnRanWu2025}.

\subsection{Spherical partial differential equations}
Hyperinterpolation has also been adapted to the numerical solution of differential equations on the sphere. Its role is to replace continuous spectral projections and inner products by discrete quantities evaluated at a spherical design. Because $\mathbb S^d$ is a closed manifold, boundary conditions are absent; depending on the differential operator, one may instead require a normalization or compatibility condition, for example, a zero-mean condition for an operator with constants in its kernel.

Let $V_n=\mathbb P_n(\mathbb S^d)$ with basis $\{p_\ell\}_{\ell=1}^{d_n}$, where $d_n=\dim V_n$, and consider
\[
   \mathscr L u=f.
\]
A discrete Galerkin formulation seeks $u_n\in V_n$ such that
\begin{equation}\label{eq:discrete-galerkin}
   \langle \mathscr L u_n,p_\ell\rangle_m
   =\langle f,p_\ell\rangle_m,
   \qquad \ell=1,\ldots,d_n,
\end{equation}
where $\langle\cdot,\cdot\rangle_m$ is evaluated by a spherical design. The design strength must be chosen so that the products occurring in \eqref{eq:discrete-galerkin} are integrated exactly, or so that the resulting quadrature error can be controlled. The required degree therefore depends on both $V_n$ and the action of $\mathscr L$; orthogonality alone does not justify the discrete system. Hyperinterpolation provides an efficient way to assemble the projected data and, in suitable settings, the projected differential operator; see \cite{wu2023restricted}.

For a time-dependent semilinear equation
\[
   \frac{\partial u(\bm x,s)}{\partial s}
   =\mathscr L u(\bm x,s)+f(\bm x,u(\bm x,s)),
\]
the nonlinear term is evaluated pointwise at the design nodes and then projected onto $V_n$ by hyperinterpolation. Thus, the nonlinearity is discretized by sampling and projection, rather than being automatically linearized. This construction leads to spectral or pseudospectral time-stepping schemes; an Allen--Cahn equation on the sphere is discussed in \cite{wu2023restricted}. The same work also studies relaxed quadrature conditions, under which exact polynomial integration is replaced by quantitative stability and consistency estimates.

\subsection{Spherical integral equations}

Hyperinterpolation has also been successfully applied to the numerical treatment of integral equations. It was initially formulated in \cite{MR878693} as a discrete orthogonal projection based on spherical $t$-designs for Fredholm integral equations of the second kind of the form
\begin{equation}
u(\bm{x}) - \int_{\Stwo} K(\bm{x},\bm{y})\, u(\bm{y})\, \dww_2(\bm{y}) = f(\bm{x}),
\label{eq:fredholm}
\end{equation}
which can be written symbolically as $(I-\Kop)u=f$. Its practical implementation relies fundamentally on spherical $t$-designs: both the quadrature nodes used to approximate the integral and the points defining the hyperinterpolation operator $\mathcal L_n$ are taken from a spherical $t$-design on $\Stwo$.

Two principal techniques, both rooted in spherical $t$-design sampling, support the use of hyperinterpolation for such equations:

\begin{enumerate}[label=\arabic*., leftmargin=*]
    \item {\em Projection-based methods.} Following the approach introduced by Atkinson and Bogomolny \cite{MR878693}, the orthogonal projection is replaced by the hyperinterpolation operator $\mathcal L_n$, constructed explicitly from a spherical $t$-design. This leads to the discrete system
    \[
    (I-\LK) u_n = \Lnf,
    \]
    where $u_n\in\Pl$ is the numerical solution. The quadrature required to evaluate $\LK$—i.e., to approximate $\int_{\Stwo} K(\bm{x},\bm{y})u(\bm{y})\,\dww_2(\bm{y})$—is performed using the same spherical $t$-design, ensuring consistency with the hyperinterpolation framework.

    \item {\em Singular-kernel methods.} For kernels that admit a factorization $K=K_1K_2$ with $K_1$ singular and $K_2$ continuous, a specialized technique can be employed. In line with the ideas of efficient hyperinterpolation, if the modified moments \eqref{equ:modifiedmoments} for the singular part $K_1$ are computable, the continuous factor $K_2u$ may be replaced by its hyperinterpolant—constructed via a spherical $t$-design—thereby avoiding direct evaluation of the singular integral. This strategy, often combined with projection-based ideas, is examined in \cite{MR2071394, MR3822242, MR1922922}. A key application concerns the singular kernel $1/|\bm{x}-\bm{y}|$ on $\Stwo$, for which Wienert \cite{Wienert1990} devised a fast algorithm that uses spherical $t$-designs as the sampling nodes, enabling efficient solutions of three-dimensional acoustic scattering problems.
\end{enumerate}

The singular-kernel approach can also be used with collocation points chosen from a spherical $t$-design, yielding a quadrature-based (Nyström) method. This is carried out in \cite{an2024spherical}, which builds on the relaxed quadrature rules of \cite{an2022quadrature} and employs spherical $t$-designs as the sampling nodes for collocation and numerical integration. A promising open direction is to incorporate relaxed-quadrature results into projection-based methods—specifically, by relaxing the polynomial-exactness requirements of spherical $t$-designs while maintaining numerical stability—thus extending the applicability of the methodology to a wider range of sampling configurations.

\subsection{Spherical data processing}
In data processing, wavelets and framelets help capture the main features of the underlying signals through sparse representations. They have been developed for Euclidean data (see, e.g., \cite{chui1992introduction, daubechies1992ten, han2017framelets, kutyniok2012shearlets, mallat1999wavelet, lin2025integral}) over the past four decades. Wavelets on the sphere first appeared in \cite{narcowich1995nonstationary, potts1995interpolatory, schroder1995spherical}. Antoine and Vandergheynst in \cite{antoine1998wavelets, antoine1999wavelets} used a group-theoretical approach to construct continuous wavelets on the sphere. Localized frames on the sphere were studied in \cite{legia2009localized, mhaskar2005representation, narcowich2006localized}. Extensions of wavelets/framelets on the sphere with more desirable properties, such as localization, tightness, symmetry, and directionality, were further studied in \cite{mhaskar2004polynomial, demanet2001directional, iglewska2017frames, wiaux2007complex, mcewen2018localisation} and many references therein.

Chen and Womersley \cite{chen2018spherical} study sparse recovery on $\mathbb S^2$ using spherical designs and nonconvex minimization. Their analysis gives error bounds for design-based sampling and replaces the $\ell_1$ penalty by an $\ell_q$ penalty with $0<q<1$ in a basis-pursuit-denoising model. The numerical study compares several spherical-design-based configurations, including $t$-design and $t_\varepsilon$-design nodes, with extremal fundamental systems and tensor-product nodes, and also compares the cases $q=1$ and $0<q<1$.

\begin{figure}[htpb!]
    \centering
    \includegraphics[width=0.32\linewidth]{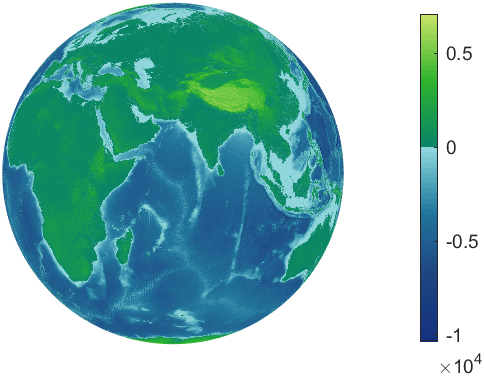}
    \includegraphics[width=0.32\linewidth]{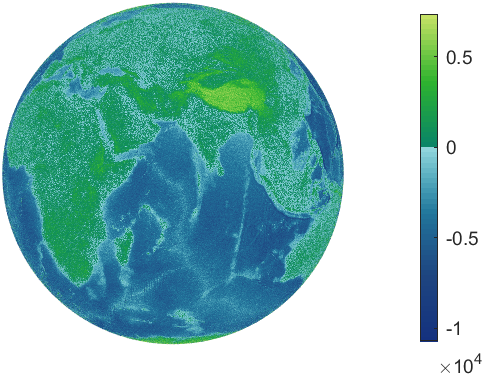}
    \includegraphics[width=0.32\linewidth]{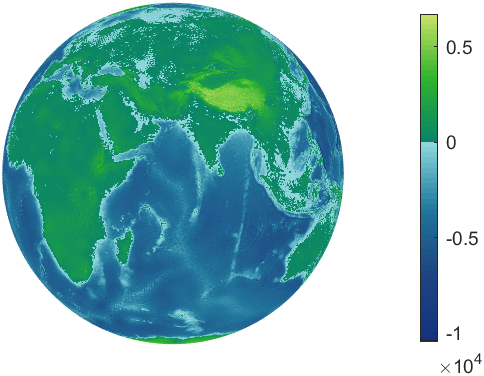}\\
    \includegraphics[width=0.32\linewidth]{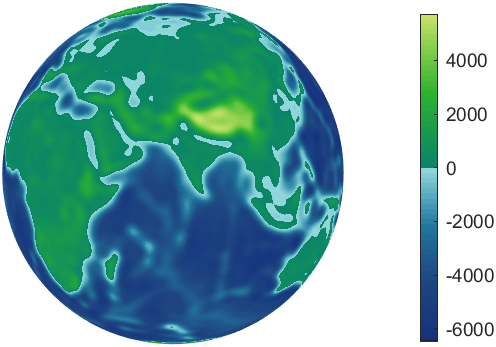}
    \includegraphics[width=0.32\linewidth]{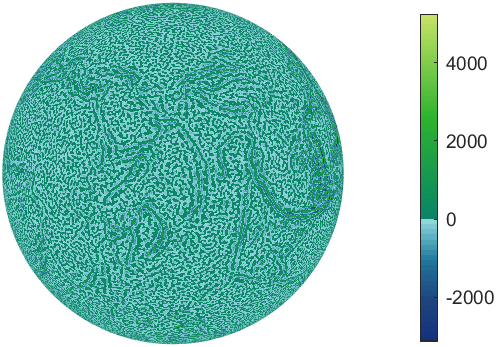}
    \includegraphics[width=0.32\linewidth]{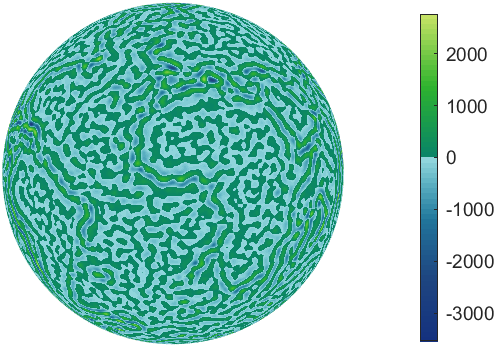}    
    \caption{ETOPO data processing. Top row (from left to right): Original, noisy, and reconstructed ETOPO data. Bottom row: decomposition of the noisy ETOPO data into low-pass filtered data (left) and high-pass filtered data at two different scales (the two panels on the right) via spherical tight framelet systems.}
    \label{fig:sph_data_processing}
\end{figure}

Using high-strength numerical spherical designs and spherical framelet systems, the authors of \cite{xiao2023spherical} perform signal and image denoising by local thresholding with a suitably chosen spherical cap restriction. Many numerical experiments are conducted to demonstrate the efficiency and effectiveness of spherical framelets and spherical designs, including Wendland function approximation, ETOPO (Earth TOPOgraphy) data \cite{amante2009etopo1} processing, and spherical image denoising. See Figure~\ref{fig:sph_data_processing} for an illustration of spherical data processing using spherical tight framelet systems.

ETOPO and CMB (cosmic microwave background) data \cite{bartlett2016planck} are prominent examples of data defined on the sphere. Spherical framelets provide localized multiscale representations for such data. They can also be combined with spherical neural-network architectures \cite{li2024convolutional}; current numerical evidence suggests useful computational possibilities, while systematic guarantees on generalization and robustness remain an open research direction. Nested spherical designs may be particularly valuable here because they permit the reuse of samples across resolution levels.

\section{Final remarks}
This survey has emphasized spherical designs as equal-weight cubature rules with exact polynomial integration, together with their algebraic and geometric structure. The modern theory contains two complementary components: existence and nonexistence results determine which cardinalities are possible, while numerical methods seek configurations that are both accurate and geometrically well distributed. Covering radius, separation, mesh ratio, and conditioning provide a common language connecting these components to approximation stability.

Several questions remain central. First, although $N_{\min}(d,t)\asymp t^d$ is known for fixed $d$, the sharp asymptotic constants and the fine structure of attainable cardinalities remain largely open. Second, certified construction close to the minimal-cardinality regime is still difficult, especially in high degrees and dimensions. Third, the optimal-order existence of nested spherical designs remains unresolved. Fourth, a quantitative theory relating polynomial exactness, covering and separation, matrix conditioning, and the stability of interpolation or hyperinterpolation would help explain when a design is not only exact but also numerically effective.

Applications in statistics \cite{haines2025optimal}, quantum information \cite{ZhuMasa2018}, differential and integral equations, and spherical data analysis continue to broaden the subject. In these settings, exactness is only one criterion: computational cost, robustness to perturbations and noise, multilevel structure, and the ability to certify a computed configuration are equally important. These considerations suggest that future work should treat spherical designs not as isolated point sets, but as structured sampling systems whose algebraic, geometric, and numerical properties must be analyzed together.

\section*{Acknowledgments}
The work of C. An was supported by the National Natural Science Foundation of China (Project No. 12371099). The work of X. Zhuang was supported in part by the Research Grants Council of Hong Kong (Project Nos. CityU 11302023 and CityU 11300825). We sincerely thank Professors  Eiichi Bannai, Peter Boyvalenkov, Renjun Duan, Yeyao Hu, Ian H. Sloan, Xiaoming Yuan, Huangjun Zhu and Ding-Xuan Zhou for their valuable comments and suggestions on the manuscript. Special thanks are extended to young scholars Yuchen Xiao, Xin Xu and Xiannan Hu for their assistance.


\end{document}